\documentclass[11pt,oneside,a4paper,mathscr]{amsart}
\usepackage{amsthm,amsmath,amssymb,euscript,amscd,a4wide}
\newtheorem{thm}{Theorem}[section]
\newtheorem{cor}[thm]{Corollary}
\newtheorem{lemma}[thm]{Lemma}
\newtheorem{prop}[thm]{Proposition}

\theoremstyle{definition}
\newtheorem{dfn}[thm]{Definition}
\newtheorem{note}[thm]{Notes}
\newtheorem{example}[thm]{Example}
\newtheorem{conjecture}[thm]{Conjecture}

\theoremstyle{remark}
\newtheorem{remark}[thm]{Remark}
\newtheorem{notation}[thm]{Notation}
{\end{note}}

\numberwithin{equation}{section}

\newcommand{\thmref}[1]{Theorem~\ref{#1}}

\newcommand{\lemref}[1]{Lemma~\ref{#1}}
\newcommand{\propref}[1]{Proposition~\ref{#1}}

\newcommand{\C}{{\mathbb C}}
\newcommand{\Z}{{\mathbb Z}}
\newcommand{\PP}{{\mathbb P}}
\newcommand{\OO}{{\mathscr O}}
\let\cross=\times
\let\tensor=\otimes

\newcommand{\Id}{\operatorname{Id}}

\newcommand{\Hilb}{\operatorname{Hilb}}

\newcommand{\Pic}{\operatorname{Pic}} 
\newcommand{\Ext}{\operatorname{Ext}}

\newcommand{\Hom}{\operatorname{Hom}}

\newcommand{\ch}{\operatorname{ch}}

\newcommand{\M}{{\mathcal M}}
\newcommand{\T}{{\mathbb T}}
\renewcommand{\P}{{\mathscr P}}
\newcommand{\FF}{{\mathcal F}}
\newcommand{\dT}{{\hat{\mathbb T}}}

\newcommand{\Ltensor}{\mathbin{\buildrel{\mathbf L}\over{\tensor}}}

\newcommand{\xhat}{{\Hat x}}

\newcommand{\I}{{\mathscr I}}

\newcommand{\refb}[2]{\cite{#1}{#2}}

\newcommand{\compo}{\raise2pt\hbox{$\scriptscriptstyle\circ$}}

\newcommand{\lra}{\longrightarrow}

\newcommand{\lRa}[1]{\>{\buildrel {#1}\over\longrightarrow}\>}
\newcommand{\lLa}[1]{\>{\buildrel {#1}\over\longleftarrow}\>}

\def\rep#1{\bysame}
\def\bib[#1]#2<#3>#4|#5(#6){\bibitem{#1}
{\sc #2},\ #3,\ {\it #4},\ {\bf #5}\ (#6)}
\def\tbib[#1]#2<#3>#4.{\bibitem{#1} {\sc #2},\ {\it #3},\ #4.}
\def\bibt[#1]#2<#3>#4|#5(#6){\bibitem{#1}
{\sc #2},\ #3,\ in\ {\it #4}\ ed.\ {#5},\ (#6)}
\def\bibit[#1]#2<#3>{\bibitem{#1} {\sc#2},\ #3}
\def\prebib[#1]#2<#3>{\bibitem{#1} {\sc #2},\ {\it #3},\ Preprint}
\def\toabib[#1]#2<#3>{\bibitem{#1} {\sc #2},\ #3,\ to appear}

\theoremstyle{definition}
\newtheorem{problem}[thm]{Problem}

\let\adj=\dashv
\newcommand{\Spl}{\operatorname{Spl}}
\newcommand{\codim}{\operatorname{codim}}
\newcommand{\EE}{\mathbb{E}}
\newcommand{\RR}{\mathbb{R}}
\newcommand{\lhom}{\operatorname{\mathcal{H}\mathnormal{om}}}
\newcommand{\R}{\mathbf{R}}
\newcommand{\natran}{\mathbin{\mathop{\longrightarrow}\limits\cmplx}}
\newcommand{\natcong}{\mathbin{\mathop{\cong}\limits\cmplx}}
\newcommand{\cmplx}{^{\scriptscriptstyle\bullet}}
\newcommand{\hcmplx}{_{\scriptscriptstyle\bullet}}
\newcommand{\Rconvol}{\mathbin{\mathop{*}\limits^{\mathbf{R}}}}
\title[]{Generalized Fourier-Mukai Transforms}
\author[]{Antony Maciocia}
\address{Department of Mathematics and Statistics\\
The University of Edinburgh\\
The King's Buildings\\ Mayfield Road\\ Edinburgh, EH9 3JZ.}
\email{ama@@maths.ed.ac.uk}
\thanks{This work was carried out with support of the Engineering and
Physical Sciences Research Council of the UK and the Consiglio Nazionale delle
Ricerche of Italy. The author is a member of the VBAC group of Europroj. The
author benefited from useful discussions with A.D.~King and U.~Bruzzo.} 
\date{\today}
\begin{document}

\begin{abstract}
We study a generalization of the Fourier-Mukai transform for
smooth projective varieties.
We find conditions under which the transform satisfies an
inversion theorem. This is done by considering a series of four
conditions on such transforms which increasingly constrain them.
We show that a necessary condition for the existence of such
transforms is that the first Chern classes must vanish and the dimensions of the
varieties must be equal. We introduce the notion of bi-universal sheaves. Some
examples are discussed and new applications are given, for example, to
prove that on polarised abelian varieties, each Hilbert scheme of points
arises as a component of the moduli space of simple bundles. The
transforms are used to prove the existence of numerical constraints on
the Chern classes of stable bundles.
\end{abstract}

\maketitle

\section*{Introduction}{}
A Fourier transform could be loosely be described as `pullback a function to
$\RR^n\cross\RR^n$ multiply by $\exp(2\pi i\langle x,y\rangle)$ and take the
direct image (integrate) with respect to the second variable'. If we regard
this as a transformation from $L^2$ integrable functions to themselves then
this satisfies certain useful properties such as having an inverse, the
Parseval theorem, the convolution theorem, etc. The Fourier-Mukai transform was
introduced in \cite{Muk1} and is formally analogous to the Fourier transform
but acts on the derived category of (bounded) complexes of sheaves on an
abelian variety $\T$ and maps this to the same category but for the dual
abelian $\dT$
variety. It too satisfies useful properties such as the Fourier Inversion
Theorem (FIT), the Parseval Theorem, the Convolution Theorem, etc (see
\cite{Muk2}). More precisely, the role of $\exp(2\pi i\langle x,y\rangle)$ is
played by the Poincar\'e line bundle over $\T\cross\dT$ and the direct image
needs to be derived.

In \cite{Muk3}, Mukai showed that a similar transform gave rise to a FIT on
the level of $K$-Theory using the universal bundle $\EE$ over $S\cross\M(S)$
instead of the Poincar\'e bundle, where $S$ is a K3 surface and $\M(S)\cong S$
is a two dimensional moduli space of simple sheaves on $S$. This was shown to
give rise to a FIT in \cite{BBH} and which satisfies the Parseval Theorem.
The question arises: when do such transformations give rise to Inversion
Theorems in more general contexts? We shall go some way to answering this
question in this paper by classifying such transformations and giving
conditions that they give rise to Inversion Theorems. We shall restrict our
attention in the examples
to the case of holomorphic varieties although it is not too hard to
extend our results to the case of varieties over more general fields.
The theorems are at their strongest in dimensions 1 and 2 but restricted forms
apply to higher dimensions as well. Our first aim is to give a rough
classification of such transforms. 

In section \ref{s:applics} we study two applications of the theory of
Fourier-Mukai transforms to the case of abelian varieties. A new
transform is constructed which is based on a certain component of
the moduli space of simple bundles on abelian varieties. This is then
used to deduce quite quickly that each Hilbert scheme of points arises
as a component of the moduli space of simple torsion-free sheaves.
This is expressed in \thmref{t:hilbismod}. We can also use these
the general theory to prove that if the torus acts freely and
effectively on any moduli component then the Euler characteristic of
the sheaves parametrised by that component must be $\pm1$ (see
\thmref{t:eulercl}).

\begin{notation}
We shall let $X$ and $Y$ be two smooth varieties and we shall consider pairs
of functors $\R\Phi:D(X)\to D(Y)$ and $\R\hat\Phi:D(Y)\to D(X)$, where $D(X)$
denotes the derived category of bounded complexes of coherent sheaves on $X$.
In this paper, for the sake of clarity, we give names to the various
possible types of such pairs. We use the terms invertible
correspondence, adjunction, Verdier and Fourier-Mukai type. These are
not mutually exclusive conditions. They are each treated in the first
four sections.

We use the notation $T_E$ to denote the functor $F\mapsto E\Ltensor F:D(X)\to
D(X)$, where $E\in D(X)$.
\end{notation}

\section{Invertible correspondences and resolutions of the diagonal}
Let $X\lLa{x}Z\lRa{y}Y$ be flat maps of smooth
quasi-projective varieties. Let $D(S)$ denote the derived category of
bounded complexes of coherent sheaves on $S$. 
Fix two objects $P$ and $Q$ in $D(Z)$ and define two functors
$\R\Phi_P:D(X)\to D(Y)$ and $\R\hat\Phi_Q:D(Y)\to D(X)$ by
$$\displaylines{\R\Phi_P(-)=\R y_*(x^*-\Ltensor P)\cr
\R\hat\Phi_Q(-)=\R x_*(y^*-\Ltensor Q).\cr}$$
Let $Z\lLa{p_y} Z_y\lRa{p'_y} Z$ (respectively $Z_x$) be the pullback of $y$
(respectively $x$) along itself. Let $q_x:Z_y\to X\cross X$ and
$q_y:Z_x\to Y\cross Y$ be the maps
defined by $(x\compo p_y,x\compo p'_y)$ and $(y\compo p_x,y\compo p'_x)$.

\begin{thm}
$\R\Phi_P$ and $\R\hat\Phi_Q$ give an equivalence
of categories (after a shift of complexes) if and only if the
following two conditions hold:\\
(i) $\R q_{x*}(p_y^* P\Ltensor p^{\prime*}_yQ)\cong
\OO_\Delta[r]$\\
(ii) $\R q_{y*}(p^{\prime*}_x P\Ltensor p_x^* Q)\cong
\OO_{\Delta'}[r]$,\\
where $\OO_\Delta$ is the structure sheaf of 
the diagonal $\Delta\subset X\cross X$,
$\OO_\Delta'$ is the structure sheaf of
diagonal $\Delta'$ in $Y\cross Y$, and $r$ is an
integer. All isomorphisms are quasi-isomorphisms of complexes.
\label{th:koszul}
\end{thm}
In other words we must have that the LHS's of (i) and (ii) are resolutions of
the both diagonals if we are to have an inversion theorem for such transforms.
The proof of the sufficiency of the conditions is, of course, well known and
fairly trivial, but the proof of necessity does require some care and so we
reproduce it here.
\begin{proof}\footnote{I am grateful to the referee for pointing out a
simplification in the original version of this proof}
Without loss of generality we set $r=0$.
Condition (i) will be equivalent to the fact that $\R\Phi_P$ has a left
inverse and (ii) will be the corresponding statement for $\R\hat\Phi_Q$.
Hence, it suffices to consider only (i). 

Now, $\R\hat\Phi_Q\compo\R\Phi_P(E)=\R x_*\bigl(y^*\R y_*(x^*E\Ltensor P)\Ltensor Q\bigr)$.
Using $Z_y$ and the base-change formula we can write this as
$$\R x_*\bigl(\R p'_{y*}(p^*_yx^* E\Ltensor p^*_yP)\Ltensor Q\bigr).$$
Using the projection formula (and the hypotheses on $X$, $Y$ and $Z$ we have
$$\R(x\compo p'_y)_*\bigl((x\compo p_y)^* E\Ltensor p^*_y P\Ltensor p^{\prime*}_y
Q\bigr).$$
Let $p_1$ and $p_2$ denote the projections $X\cross X\to X$. Then $x\compo
p_y=p_1\compo q_x$ and $x\compo p'_y=p_2\compo q_x$. Substituting these and
using the projection formula again we have
\begin{equation}
\R p_{2*}\bigl(p_1^*E\Ltensor\R q_{x*}(p^*_y P\Ltensor p^{\prime *}_y
Q)\bigr).\label{eq:psiophi}
\end{equation}
Suppose condition (i) holds then
$$\R\hat\Phi_Q\compo\R\Phi_P=\R p_{2*}(p^*_2 E\Ltensor\OO_{\Delta})=E.$$
This proves the sufficiency of (i). To see
that it is necessary observe that \ref{eq:psiophi} still holds. Assume that
$\R\hat\Phi_Q\compo\R\Phi_P\natcong\Id$ and put $E=\OO_\alpha$, the structure
sheaf of a point $\alpha\in X$. Let $A\hcmplx$ be a (bounded) complex of
locally-free $\OO_{X\cross X}$-modules representing
$\Gamma=\R q_{x*}(p^*_y P\Ltensor p^{\prime *}_y Q)$.
By assumption, we have the quasi-isomorphism
\begin{equation}
\R p_{2*}(\Gamma\Ltensor p_1^*\OO_\alpha) \simeq\OO_\alpha.\label{eq:alpha}
\end{equation}
Then
$$R^*p_{2*}(\Gamma\Ltensor
p_1^*\OO_\alpha)=H^*\bigl(p_{2*}(A\hcmplx\tensor p_1^*\OO_\alpha)\bigr).$$
But this is concentrated in position 0 and so we
see that $H_i(A\hcmplx\tensor p_1^*\OO_\alpha)=0$  for all $i\neq0$ and so
$\Gamma$ is flat over $p_1$.
Now the fibres of $\Gamma$ are given by
$\Gamma_{(\alpha,\beta)}=
\R p_{2*}(\Gamma\tensor p_1^*\OO_\alpha\tensor p_2^*\OO_\beta)$. We
can rewrite this as $\R p_{2*}(\Gamma\tensor p_1^*\OO_\alpha)\tensor\OO_\beta
=\OO_\alpha\tensor\OO_\beta$. This is zero
if $\alpha\neq\beta$. Hence, $\Gamma$ is supported on
$\Delta_X$ and the case $\alpha=\beta$ implies that it has rank 1
everywhere along this diagonal. If we now substitute $E=\OO_X$ then a standard
hypercohomology argument shows that it must also be trivial.
\end{proof}
\begin{dfn}
We shall call functors $\R\Phi$ and $\R\hat\Phi$ a {\it invertible
correspondence transforms}
if they satisfy the conditions of \thmref{th:koszul}. These can
also be thought of as transformations satisfying the `Fourier Inversion Theorem'
or FIT.
\end{dfn}
The invertible correspondence transforms should be compared to the notion
of tilting transforms for categories of modules over associative
rings. Details of this can be found in \cite{Ric}.

As an intermediate example one can consider the Beilinson spectral sequence
which could be viewed as a  composition of derived functors from the
derived category of coherent sheaves on $\C P^n$ to the derived
category of finitely generated modules over a suitable algebra (the
path algebra of a certain quiver) see \cite{Bei1} and \cite{Bei2}. This also
works for more general varieties, see King \cite{Ki}.

We shall see in section 3 that when we consider invertible
correspondences for varieties then their existence imposes quite
strong condition on the varieties. 

\section{Transforms of Adjunction Type}
In this section we shall look briefly at the categorical aspects of
equivalences of such transforms. The following is well known:
\begin{thm} {\rm(See \cite{McL})}. Let $F:\mathcal{A}\to\mathcal{B}$ and
$G:\mathcal{B}\to\mathcal{A}$ be two functors of categories. If
$F$ and $G$ give rise to an equivalence of categories (i.e. $F\compo G$
and $G\compo F$ are naturally equivalent to their respective identity
functors) then $F\adj G\adj F$. Conversely, if $F\adj G\adj F$ and $F$ is fully
faithful such that it surjects on objects up to isomorphism (we shall say {\em
quasi-surjects}) then $F$ and $G$ determine an equivalence of categories.
\end{thm}
\begin{dfn}
We say that a pair of functors $F$ and $G$ are a 
{\it transform of adjunction type} if $F\adj G\adj F$.
\end{dfn}
The theorem says that invertible correspondence transforms are adjunction
transforms. 
When we have an adjunction type transform then we often already know that the
functors are inverses on one side. For example, suppose that we already know
that the adjunction gives rise to
$GF\mathrel{\mathop{\cong}\limits\cmplx}\Id$. Recall that an adjunction gives
rise to natural transformations $\eta:\Id\natran GF$ and
$\epsilon:FG\natran\Id$ called the unit and counit of the adjunction $F\adj G$
similarly $G\adj F$ gives rise to $\eta':\Id\natran FG$ and
$\epsilon':GF\natran\Id$. Then it is well known that $G$ is faithful if and
only if $\epsilon_a$ is epi for all $a$ and $G$ is full if and only if
$\epsilon_a$ is split monic (i.e. has a left inverse). The same holds for $F$
if we replace $\epsilon$ by $\epsilon'$. If we already know that $\epsilon'$
is a natural isomorphism then $F$ must be full and faithful and that $G$ is
full. But $GFa\cong a$ for each $a$ and so $G$ must quasi-surject on objects.
It now follows that $F$ and $G$ are an equivalence of categories if and only
if $G$ is faithful. We summarise this in the following.
\begin{prop}\label{p:parseval}
Suppose that $F$ and $G$ form a transform of adjunction type such that
$GF\natran\Id$ is an isomorphism. Then $F$ satisfies the {\em Parseval
Theorem}: $$\Hom(a,b)\cong\Hom(Fa,Fb)$$
Furthermore, $F$ and $G$ form an equivalence of
categories if and only if $G$ satisfies the Parseval Theorem.\label{p:semiadj}
\end{prop}

\section{Transforms of Verdier Type}
We now look at more specific transforms. We limit ourselves to smooth
projective varieties $X$ and $Y$ and let $Z=X\cross Y$.
Let $n=\dim X$ and $m=\dim Y$. Note that we have 
$$\R\hat\Phi(F)=\R x_*\R\lhom(P,y^*F),$$
where $Q=\R\lhom(P,\OO)$ and $P$ is a sheaf. We shall see that this choice of
$Q$ is forced on us up to the pullback of (a power of) the canonical bundle
and a shift.
Recall that we have Grothendieck-Verdier duality:
$$\Hom_{D(X)}(\R x_*F,G)\cong\Hom_{D(X\cross Y)}(F,x^*G\Ltensor
y^*\omega_Y[n]).$$ 
$$\Hom_{D(Y)}(\R y_*F,G)\cong\Hom_{D(X\cross Y)}(F,y^*G\Ltensor
x^*\omega_X[m]).$$
(see Hartshorne \refb{Hartres}{III.11}).
Applying this and the classical adjunction $f^*\adj\R f_*$ to $\R\Phi$ and
$\R\hat\Phi$ we obtain
\begin{prop}
$$\R\Phi\adj T_{\omega_X}\compo\R\hat\Phi[n]
\quad\text{and}\quad
\R\hat\Phi\adj T_{\omega_Y}\compo\R\Phi[m],$$
where $Q=\R\lhom(P,\OO_{X\cross Y})$. 
\end{prop}
\begin{proof}
This is just a computation using the adjunctions above:
\begin{align*}
\Hom_{D(Y)}(\R y_*(x^*F\Ltensor P),G)&\cong
\Hom_{D(X\cross Y)}(x^*F\Ltensor P,y^*G\tensor x^*\omega_X[n])\\
&=\Hom_{D(X)}(F,\R x_*\R\lhom(P,y^*G)\tensor\omega_X[n])
\end{align*}
Similarly for the other adjunction.
\end{proof}
\begin{remark}
We could use $Q=\R\lhom(P,x^*L\tensor y^*M)$ for line bundles $L$ on $X$ and
$M$ on $Y$. This does not affect the adjunctions except we have to twist by
one or other of these line bundles or their duals.\label{r:libQ}
\end{remark}
\begin{prop}\label{p:veradj}
Let $\R\Delta_X=\R\lhom(-,\omega_X)[n]$ and
$\R\Delta_Y=\R\lhom(-,\omega_Y)[m]$. Then
$$\R\Delta_Y\compo\R\Phi_P\natcong T_{\omega_Y}\compo
\R\Phi_{Q}\compo\R\Delta_X[m],$$
$$\R\Delta_X\compo\R\hat\Phi_P\natcong T_{\omega_X}
\compo\R\Hat\Phi_{Q}\compo\R\Delta_Y[n],$$
where $Q=\R\lhom(P,\OO)$
\label{p:dual}
\end{prop}
\begin{proof}
We use local Verdier duality:
\begin{align*}
\R\lhom(\R y_*(x^*F\Ltensor P),\omega_Y)&\cong
\R y_*\R\lhom(x^*F\tensor x^*\omega^*_X,\R\lhom(P,y^*\omega_Y[n+m]))\\
&\cong\R y_*(x^*\R\lhom(F,\omega_X)\Ltensor\R\lhom(P,y^*\omega_Y)[n+m])\\
&\cong\R y_*(x^*\R\lhom(F,\omega_X[n])\Ltensor\R\lhom(P,\OO))[m]\tensor
\omega_Y
\end{align*}
Similarly for $\R\hat\Phi$.\end{proof}
\begin{thm}\label{c:koko}
Let $X$ and $Y$ be two smooth projective varieties of dimensions $n$
and $m$ respectively. Let $P$ and $Q$
be two complexes of coherent sheaves on $X\cross Y$. We assume that $X\neq Y$
and $P$ is a sheaf of rank at least 1.
Define two functors $D(X)\to D(Y)$ and $D(Y)\to D(X)$ by
$\R\Phi(E)=\R y_*(x^*E\Ltensor P)$ and $\R\hat\Phi(F)=\R
x_*(y^*F\Ltensor Q)$ respectively, where $x$ and $y$ are the
projections from $X\cross Y$ to $X$ and $Y$ respectively.
If $\R\Phi$ and $\R\hat\Phi$ form an equivalence of categories then 
$\omega_X^k=\OO_X$, $\omega_Y^k=\OO_Y$ for some integer $k$, $n=m$ and
$Q=\R\lhom(P,x^*\omega_X)[n]=\R\lhom(P,y^*\omega_Y)[n]$. In particular, we
must have $c_1(X)=0$.
\end{thm}
\begin{proof}
This follows from \propref{p:veradj} and the
fact that adjoints are unique up to natural isomorphism. We may assume
that $Q$ takes the form $\R\lhom(P,x^*\omega_x)[n]$ because the adjoint of $\R\Phi_P$
is $\R\Hat\Phi_{\R\lhom(P,x^*\omega_X)[n]}$ by duality.
The adjoints give
$$\R\lhom\bigl(\R\lhom(P,\OO),\OO\bigr)\tensor x^*\omega^*_X\tensor
y^*\omega_Y[m-n]=P$$ 
which can only happen if $n=m$ and, since $P$ has support on the whole of the product,
$\omega_X^k=\OO_X$ and $\omega_Y^k=\OO_Y$ for some integer $k$ as
required. In fact, what we need is that $P\tensor x^*\omega_X\tensor y^*\omega^*_Y\cong
P$.
\end{proof}
\begin{remark}
In other words, an invertible correspondence transform given by, for example,
a torsion-free sheaf $P$ for smooth
projective varieties can only exist if the dimensions of the varieties
are the same and their canonical bundles are trivial or torsion (with order
dividing the rank of $P$). Moreover, the
transform complex $Q$ must be the (derived) dual of $P$ (up to a twist and shift).
Relative versions of the results of this section
are also available where we we have $Z=X\cross_S Y$ for some scheme $S$. In
fact, we can clearly strengthen the conclusion to say that the canonical
bundles must act trivially on the the restrictions of $P$ to
$X\cross\{y\}$ and $\{x\}\cross Y$. This will happen, for example, for
relative transforms where the fibres have trivial canonical bundle.
\end{remark}
\begin{dfn}
We say that a pair of transforms $\R\Phi_P$ and $\R\hat\Phi_Q$ is {\em of
Verdier type} if $\omega_X=\OO_X$, $\omega_Y=\OO_Y$, $\dim X=\dim
Y$ and $Q=\R\lhom(P,\OO_{X\cross Y})$ (up to a shift).
So we have the chain of implications 
$$\hbox{invertible correspondence}\implies\hbox{Verdier}
\implies\hbox{adjunction}$$.
\end{dfn}
\section{Fourier-Mukai Transforms}
\begin{dfn} Let $\R\Phi$ and $\R\hat\Phi$ be transforms of Verdier
type as in the last section but with the additional constraint that
$P$ is a locally-free sheaf over $X\cross Y$.
We also assume that they give rise to an equivalence of categories.
Then we say that
they are {\it transforms of Fourier-Mukai type}. In other words, the
transforms are both of invertible correspondence type and Verdier type
with some constraint on the choice of $P$.
\end{dfn}
Let $\Spl(S)$ denote the moduli space of simple sheaves on $S$. Recall that
Mukai \cite{Muk3} proves that this space is smooth when $S$ is an abelian
surface or a K3 surface.
\begin{prop}\label{c:nisr}
Suppose that $\R\Phi$ 
and $\R\hat\Phi$ form a pair of functors which are of Fourier-Mukai
type. Then $n=r$, where $n=\dim X$ and $r$ is given in
\thmref{th:koszul}. Furthermore,
if $X$ and $Y$ admit smooth moduli of simple torsion-free sheaves then
$X$ is a component of $\Spl(Y)$ and $Y$ is a component of $\Spl(X)$.
\end{prop}
\begin{proof} In the proof we only assume that $P$ is torsion-free and flat
over both projections.
For this we use Proposition 2.26 of \cite{Muk4}. This states that if $f:Z\to
Y$ is a proper map of noetherian schemes and $F$ is a coherent sheaf on $Z$
flat over $Y$ and $S\subset Y$ is a locally complete intersection then if
$H^i(f^{-1}(y),F_y)=0$ for all $i<\codim S$ and $y\not\in S$ then $R^if_*F=0$
for the same set of $i$. 
We apply this to $S=\Delta\subset X\cross X$. Then
$\codim S=n$ and so if $r<n$ we have that the cohomology sheaves
$R^iq_{x*}(\Gamma)$ of $\R q_{x*}(\Gamma)$ vanish for for $i>r$,
where 
$$\Gamma=\R\lhom(p^*_xP,p^*_yP).$$ 
Hence $(R^rq_{x*}(\Gamma))_{(a,a')}\cong H^r(Y;\Gamma_{(a,a')})=0$ for
$(a,a')\in X\cross X\setminus\Delta$. Note that
$\Gamma_{(a,a')}\cong\R\lhom(P_a,P_{a'})$. 
This implies that $R^rq_{x*}(\Gamma)=0$
and so $R^iq_{x*}(\Gamma)=0$ for all $i$; a contradiction. Hence
$n=r$. Let $P_a$ denote the restriction of $P$ to $\{a\}\cross Y$.
We know that $h^n(Y;\lhom(P_a,P_a))=\dim\Ext^n(P_a^{**},P_a)=1$
and so $\dim\Ext^n(P_a,P_a)=1$ because $P_a\to P_a^{**}$ induces a
surjection on Ext groups. Then Serre duality implies that $P_a$ is simple.
Similarly $P_b$ is simple for all $b\in Y$. On the other hand,
$\Ext^n(P_a,P_{a'})=0$ for $a\neq a'$ and so $P_a\not\cong P_{a'}$ and hence
$X\subset\Spl(Y)$. 

Suppose that $E\in\Spl(Y)$ is not in $X$ but has the same Chern character as
$P_a$. If $\Ext^i(P_a,E)=0$ for all $i$
and all $a\in X$ then $\R\hat\Phi(E)=0$. But $\R\hat\Phi$ is an equivalence of
categories and so this is impossible. For a fixed $a$ the set
$Z_a=\{E\in\Spl(Y):\Ext^i(P_a,E)=0,\;\forall i\}$ is Zariski open in $\Spl(Y)$
and non-empty as it contains $X\setminus\{a\}$. But $(U_a\cap Z_a)\setminus X$
is empty for any open neighbourhood $U_a$ of $P_a$ in $\Spl(Y)$. This shows
that a component of $\Spl(Y)$ which contains $X$ is smooth at each point of
$X$ and has the same dimension as $X$.
\end{proof}
\begin{remark}
In the case when $n=1$ or $n=2$ we can proceed more directly and explicitly.
Assume for simplicity that $P$ is locally-free and consider first the case
$n=1$. 
Note that $X$ must be an elliptic curve by \ref{c:koko}.
The fact that $\R\Phi$ and $\R\hat\Phi$ are invertible correspondences
mean that 
$$\R q_{x*}\R\lhom(p^{\prime*}_y P^*,p_y^* P)=\R
q_{x*}\Lambda\cong\OO_\Delta[-r]$$ 
for some $r=0$ or 1.  Serre duality implies that $R^1 q_{x*}\Lambda\neq0$ and
so $r=1$. Then $\chi(P_b,P_c)=0$ for $c\neq b$ and hence for $b=c$ as well. We
also have $\dim\Ext^1(P_b,P_b)=1$ and so $P_b$ are simple sheaves with a 1
dimensional moduli which contains $Y$. Hence, $\dim Y=1$ as required.

When $n=2$ we argue similarly.
Then again $\chi(P^*_a\tensor P_{a'})=0$ for all $a,a'\in X$.
Note that \ref{c:koko} shows that the canonical bundles are trivial and so
Serre duality implies that $r\neq0$; otherwise $H^0(P_b^*\tensor P_b)\neq0$
so $H^2(P_b^*\tensor P_b)\neq0$, and hence $R^2q_{x*}\Lambda\neq0$. But if
$r=1$ then the support of $\{b\in Y:H^1(X,P^*_b\tensor P_c)\neq0\}$
is 0-dimensional Mukai's result above implies that $R^1q_{x*}\Lambda=0$, a
contradiction. Hence we must have $r=2$. It also follows that $P_b$ are all
simple. Then $Y$ is contained in the moduli space of such $P_b$'s and
hence has dimension at most 2 because $\chi(P_b,P_b)=0$. It cannot have
dimension 1.
\end{remark}

\section{Bi-universal sheaves}{}
In this section we ask the following question. If $X$ is a smooth complex
projective variety with trivial canonical bundle and
$Y=\M(X)$ is a smooth moduli of simple sheaves on $X$ (also assumed to have
trivial canonical bundle) then when do we obtain a
Fourier-Mukai transform from $X$ to $Y$?
We answer this question by considering various universal sheaves on the
product $X\cross Y$.
\begin{dfn}
Following Mukai, we say that a sheaf $\EE$ is {\it semi-universal} on $X\cross
Y$ if for all $y\in Y$, $\EE_b\cong E^{\oplus\sigma}$
for some $\sigma\in\Z$ and $\EE$
is a universal deformation. If $\sigma=1$ we say that $\EE$ is a {\it universal
sheaf} (as usual). If $X$ is isomorphic to a moduli of simple sheaves on $Y$
such that $\EE_a\cong E^{\oplus\sigma'}$ and $\EE_b\cong E^{\oplus\sigma}$
we say that $\EE$
{\it is bi-semi-universal}. If either $\sigma$ or $\sigma'$ is 1 then we call
$\EE$ {\it sesqui-universal} and if $\sigma'=\sigma=1$ then we call $\EE$ {\it
bi-universal}. In all these cases we say that $\EE$ is {\em strongly}
universal (etc.) if whenever $b\neq b'$ then
$\Ext^i(\EE_b,\EE_{b'})=0$ for all $i$ and similarly for $a\neq a'$.
\end{dfn}
Mukai has shown that if $Y$ is a (representable) component of the moduli of
simple sheaves on
$X$ then a semi-universal sheaf always exists (see
\refb{Muk3}{Thm~A.5}). 
\begin{remark} Of course, in dimension 2, if the component of the moduli
space consists of, say, stable bundles then the strong condition will always
hold. 
\end{remark}
In the following we let $P=\EE$. We also assume that $\dim Y=\dim X=n$.
\begin{prop} If $\EE$ is strongly semi-universal then
$\R\Phi\compo\R\hat\Phi\natcong\Id^{\oplus\sigma^2}[n]$.
In particular, $\R\hat\Phi$ is faithful. If $\EE$ is strongly
bi-semi-universal then $\sigma=\sigma'$.
\end{prop}
\begin{proof}
The conditions on $\EE$ ensures that $R^iq_{x*}\Gamma=0$ for $i<n$ and for
$i=n$ is supported on the diagonal. That it is trivial and given by such a
direct sum follows from the fact that
$\R\hat\Phi(\OO_b)\cong(\EE_b^{\oplus\sigma})^*[n]$ and
$\R\Phi(\EE_b)=\OO_b^{\oplus\sigma}$. Then $\R\Phi\compo\R\hat\Phi$ is
faithful and so $\R\hat\Phi$ is.
If $\EE$ is bi-universal then we also have
$\R\Phi(\EE^*_b)\cong\OO_b^{\oplus\sigma}$. So
$\R\Phi\compo\R\hat\Phi(\OO_b)=\R\Phi(\EE_b^{\oplus\sigma})^*=
\OO_b^{\oplus\sigma^2}$
On the other hand, $\R\Phi\compo\R\hat\Phi=\Id^{\oplus\sigma'{}^2}$ and so
$\sigma^2=\sigma'{}^2$. 
\end{proof}
\begin{cor}
If $\EE$ is strongly sesqui-universal then it must be strongly bi-universal.
\end{cor}
\begin{cor}
From \ref{p:semiadj} we see that if $\EE$ is strongly universal then the
following are equivalent
\begin{enumerate}
\item $\EE$ gives rise to a Fourier-Mukai transform.
\item $\EE$ is sesqui-universal.
\item $\R\Phi$ satisfies the Parseval Theorem. 
\end{enumerate}
\end{cor}

\section{Examples}{}
\begin{example}
We shall now look at some examples of Fourier-Mukai transforms.
The first is the Mukai transform itself. For this we set $P=\P$, the
Poincar\'e bundle on $\T\cross\dT$, where $X=\T$ is an abelian variety and
$Y=\dT\cong\Pic^0\T$ is its dual abelian variety. It was shown in \cite{Muk2}
that $\R\FF=\R\Phi_\P$ and $\R\hat\FF=\R\hat\Phi_\P$ are of
Fourier-Mukai type. In
this particular case we also obtain a convolution theorem (see \refb{Muk2}{})
as well. Note also that $X\neq Y$. A relative version can also be found in
\cite{Muk4}. The Grothendieck-Riemann-Roch Theorem can be used to compute the
Chern characters of the transforms:
$$\ch(\R\FF(E))_i=(-1)^i\ch(E)_{n-i},$$
where we identify $H^i(\T)$ with $H^{n-i}(\dT)$ via Poincar\'e
duality. It is conventional in this case to use $\R\hat\Phi_{\P}$
instead of $\R\hat\Phi_{\P^*}$.
\end{example}
\begin{example}
The second example first appeared in \cite{Muk4} and was shown to be a
Fourier-Mukai transform by Bartocci, Bruzzo and Hern\'andez Ruip\'erez
\cite{BBH}. In this case $X$ is a K3 surface satisfying suitable conditions
for the existence of $Y$, a 2-dimensional space of stable bundles on $X$ which
is isomorphic to $X$. Then $P$ is a bi-universal bundle on $X\cross Y$. The
Grothendieck-Riemann-Roch Theorem also tells us the Chern characters.
\end{example}
\begin{example}
Consider a smooth projective variety $X$ with trivial canonical
bundle and $h^{0,p}=0$ for $p\neq0,n$ (for example, a K3 surface).
Let $P=\I_\Delta$ be the ideal sheaf of the diagonal in
$X\cross X$. Then $P$ is strongly bi-universal. This follows because
$P_a\cong\I_a$ and, if $a\neq a'$ then $\Ext^i(\I_a,\I_{a'})=0$ for
all $i$ as can be easily seen using the long exact sequences induced
by the structure sequences of $a$ and $a'$.
Note that $Q$ cannot be represented by a single sheaf. The Chern
character of the resulting Fourier-Mukai transform of $E\in D(X)$ is
$$\bigl(\chi(E)-\ch_0(E),-\ch_{1}(E),\ldots,-\ch_{i}(E),\ldots,-\ch_{n-1}(E),-\ch_n(E)\bigr),
$$
as can be easily seen from the Grothendieck-Riemann-Roch formula or directly
from the structure sequence of $\Delta$.
\end{example}
\begin{example}\label{e:Mrl1}
Consider an abelian surface $\T$ polarised by $\ell$ such that $\ell^2=2r$,
with dual polarisation $\hat\ell$ of $\dT$.
Let $\M=\M(r,\ell,1)$ be the moduli space of stable bundles with the given Chern
characters. Such module spaces have been also considered by Mukai
(\cite{Muk1}). We shall prove in \propref{p:7.1} below that this is
projective and non-empty. It is easy to see that if $E\in\M$ then $R\FF(E)=
\hat L^*\tensor\P_x$ for
some (assumed symmetric) $\hat L\in\hat\ell$ and $\P_x\in\Pic\dT\cong\T$.
This implies
that $\M\cong\T$ under $E\mapsto x$. A result of Mukai (\refb{Muk3}{Appendix
2}) implies that a universal sheaf $\EE$ exists over $\T\cross\M$. In fact it
is possible to write this down explicitly. Let $\pi_{ij}$ be the projection
maps to the $i$th and $j$th factors of $\T\cross\dT\cross\T$ and $\pi_i$, the
projections to the $i$th factor. Then it is easy
to check that
$\EE=\R\pi_{13*}(\pi_2^*\hat L^*\tensor\pi_{23}^*\P\tensor\pi_{12}^*\P^*)$.
This also 
shows that $\EE$ is bi-universal. Then $\R\Phi_{\EE}$ is of
Fourier-Mukai type. Using Grothendieck-Riemann-Roch
(or \lemref{l:phiisf}) we find
\begin{align*}
\ch(\R\Phi(E))_0&=\ch_0(E)+\ch_1(E)\cdot\ell+r\ch_2(E),\\
\ch(\R\Phi(E))_1&=\ch_1(E)+\ch_2(E)\ell,\\
\ch(\R\Phi(E))_2&=\ch_2(E).
\end{align*}
We shall study this example in more detail below.
\end{example}
\section{Two Applications}\label{s:applics}
We shall now look at two applications of the general theory of Fourier-Mukai
transforms. The first is a special case of Example \ref{e:Mrl1} and
the second is a generalisation of that example.

First, we let $X=\T$ be an abelian surface with a polarisation $\ell$.
Let $\M=\M(r,\ell,1)$ denote the moduli space of (Gieseker) stable
bundles of the given Chern character. Choose symmetric representative line
bundles $L\in\ell$ and $\hat L\in\hat\ell$ in $\ell$ and the dual
polarisation respectively.
\begin{prop}\label{p:7.1}
The moduli space $\M$ is isomorphic to $\T$.
\end{prop}
\begin{proof}
Consider the collection $\{\R\hat\FF(\hat L\tensor\P_x)\mid x\in\T\}$.
Note that $R^i\hat\FF(\hat L)=0$ unless $i=0$ and so this set consists of
vector 
bundles of Chern character $(r,\ell,1)$. Moreover, these are all
$\mu$-stable since the projective bundle $B=\PP R^0\hat\FF(\hat L)$ has fibres
consisting of the linear systems $|\hat L\tensor\P_x|$ which are just
translates of $|\hat L|$. So $B$ admits a flat connection
given by translation of this linear system. This connection then
induces a projectively anti-self-dual connection on each of the
elements of the collection. This implies that the bundles are all
$\mu$-stable. Since the collection is non-empty and the Mukai
transform is an isomorphism of schemes we see that $\M\cong\T$.
\end{proof}
We have seen in \ref{e:Mrl1}
that a strongly bi-universal sheaf $\EE$ exists over $\T\cross\M$. We shall
prove the following 
\begin{thm}\label{t:modideal}
Let $(\T,\ell)$ be a polarised torus with $\ell^2=2r$.
There is a component of the moduli space of simple sheaves over $\T$ with Chern characters
$(rn-1,n\ell,n)$ is canonically isomorphic to $\Hilb^n\T\cross\dT$.
\end{thm}
We shall see that this isomorphism is given by $\R\Phi$.
We introduce the following terminology, again following Mukai.
\begin{dfn}
We say that a sheaf $E$ on $X$ satisfies $\Phi$-WIT$_i$ if for all $j\neq i$,
$R^j\Phi(E)=0$. In other words, $\R\Phi(E)$ is again a sheaf. We just write
WIT$_i$ for $\FF$-WIT$_i$ and $\Phi$-WIT if we don't want to specify $i$.
\end{dfn}
\begin{lemma}\label{l:PhiofP}
The flat line bundles $\P_\xhat$ satisfy $\Phi$-WIT$_0$ and
$R^0\Phi(\P_\xhat)=\P_\xhat$. 
\end{lemma}
\begin{proof}
This follows from Lemma \ref{l:phiisf} below since $\P_\xhat$ satisfies
WIT$_2$ with transform $\OO_{-\xhat}$.
\end{proof}
\begin{lemma}
The ideal sheaf $\I_S$ of a 0-dimensional subscheme $S$ of $\T$ satisfies
$\Phi$-WIT$_1$ and the transform can be written as $A/\OO$, where $A$ admits a
filtration whose factors are elements of $\M(r,\ell,1)$. More generally,
$\R\Phi(\I_S\tensor\P_\xhat)=A/\P_\xhat$.
\end{lemma}
Note that $\I_S$ never satisfies WIT. We could say that $\I_S$ is a
``half-WIT'' sheaf. 
\begin{proof}
Observe first that the structure sheaf $\OO_S$ of $S$ satisfies $\Phi$-WIT$_0$
and its transform is a sheaf $A$. Since $\OO_S$ is built up from a series of
extensions of structure sheaves of single points we see that $A$ admits a
filtration whose factors are $\R\Phi(\OO_x)=E_x$ for $x\in S$.
If we apply $\R\Phi$ to the structure
sequence of $S$ twisted by $\P_\xhat$ then we obtain the long exact sequence
$$0\lra R^0\Phi(\I_S\tensor\P_\xhat)\lra\P_\xhat\lRa{f} A\lra
R^1\Phi(\I_S\tensor\P_\xhat)\lra0$$
using \lemref{l:PhiofP}. Since $f=\R\Phi(\P_\xhat\to\OO_S)$ it must be non-zero
as $\P_\xhat$ and $\OO_S$ both satisfy $\Phi$-WIT. Since $A$ is locally-free
and the rank of $\P_\xhat$ is one we see that $f$ must inject. This proves the
lemma.
\end{proof}
Observe that $\ch(R^1\Phi(\I_S))=(r|S|-1,|S|\ell,|S|)$. To prove the theorem
it suffices to show that $R^1\Phi(\I_S\tensor\P_\xhat)$ is simple. But this
follows immediately from the Parseval Theorem (\ref{p:parseval})
since $\I_S$ is simple. Then the map
$\I_S\tensor\P_\xhat\mapsto\R\Phi(\I_S\tensor\P_\xhat)[1]$ gives an
injection $\Hilb^n\T\cross\dT\to\Spl(rn-1,n\ell,n)$. Since the moduli
of simple sheaves on $\T$ with this Chern character is smooth of
dimension $2n-2$ then the image must be single reducible component.
Since $\R\Phi$ is an
equivalence of derived categories, it also preserves the holomorphic
deformation structure of the spaces. In particular, the tangent spaces are
canonically isomorphic and so the map is a diffeomorphism and preserves the
complex structures. One can also see this algebraically by observing that the
Fourier-Mukai transforms are actually a natural isomorphism of moduli functors
and so give an isomorphism of coarse moduli schemes.
In fact, the transform also preserves the natural symplectic
structures which are simply given as a composition of derived morphisms. This
completes the proof of the theorem.$\qquad\square$

\medskip
For our second application, we consider an abelian variety $\T$ of any
dimension $n$ and consider a moduli space $\M$ of stable bundles on $\T$
of dimension $n$ which is isomorphic to $\T$. Suppose further that
this isomorphism is given by the translation action of $\T$ on $\M$ by
pullback. Then Mukai has shown that there is a (strongly) semi-universal sheaf
on $\T\cross\M$. We let $E_0$ be the base point determined by
$0\in\T\cong\M$. 
\begin{prop}
Define a complex $\EE$ in $D(\T\cross\T)$ given by
$$\EE=\R\pi_{13*}(\pi_2^*E_0^*\tensor\pi_{23}^*\P^*\tensor\pi_{12}^*\P),$$
where $\pi_{ij}$ denotes the projection
maps to the $i$th and $j$th factors of $\T\cross\dT\cross\T$,
$\pi_i$ is the projection to the $i$th factor and $\P$ is the
Poincar\'e bundle on $\T\cross\dT$. If we assume that $E_0$ satisfies
WIT then $\EE$ is a bi-universal sheaf.
\end{prop}
This proposition follows more or less immediately from the following
lemma.
\begin{lemma}\label{l:phiisf}
Let $\EE$ be given as in the proposition (we do not assume that $E_0$
satisfies WIT) then
$$\R\Phi\natcong(-1)^*\R\hat\FF\compo T_{\R\FF(E_0)}\compo\R\FF$$
and
$$\R\hat\Phi\natcong(-1)^*\compo\R\hat\FF\compo
T_{\R\FF(E^*_0)}\compo\R\FF[n]$$
\end{lemma}
\begin{proof}
These are just hideous computations of derived functors. Let $\hat
E=\R\FF(E_0)$. Consider the commuting diagram
\begin{equation*}
\begin{CD}
\T\cross\dT @<\pi_{12}<<  \T\cross\dT\cross\T  @>\pi_{23}>> \dT\cross\T\\
@VpVV   @VV\pi_{13}V    @VV p V \\
\T  @<x<<  \T\cross\T   @>y>> \T
\end{CD}
\end{equation*}
Then 
\begin{equation}\label{eq:aa}
\R y_*(x^*F\Ltensor\EE)=
\R y_*(x^*F\Ltensor\R\pi_{13*}
(\pi_2^*\hat E\tensor\pi^*_{23}\P^*\tensor\pi_{12}^*\P)
\end{equation}
Note that $x^*F=x^*\R p_*p^*F=\R\pi_{13*}\pi^*_{12}p^*F=\R\pi_{13*}\pi_1^*F$
so 
$$\hbox{(\ref{eq:aa})}=\R y_*\R\pi_{13*}(\pi_1^*F\tensor\pi_2^*\hat E
\tensor\pi_{23}^*\P^*\tensor\pi_{12}^*\P).$$
Which we can write as 
$$\R
p_*\bigl(\R\pi_{23*}(\pi_1^*F\tensor\pi_{12}^*\P)\tensor\P^*\tensor 
q^*\hat E\bigr)$$
since $\pi_2=q\compo\pi_{23}$.
But $\R\pi_{23*}(\pi_1^* F\tensor\pi_{12}^*\P)=
\R\pi_{23*}\pi_{12}^*(p^*F\tensor\P)=q^*\R q_*(p^*F\tensor\P)
=q^*\R\FF(F)$. Then
$$\hbox{(\ref{eq:aa})}=\R p_*\bigl(q^*(\R\FF(F)\tensor\hat E)\tensor\P^*\bigr)$$
as required.

The second equation follows similarly if we observe that
$$\R\lhom(\EE,\OO)=\R\pi_{13*}(\pi_2^*\R\lhom(\hat E,\OO)\tensor\pi_{23}^*\P
\tensor\pi_{12}^*\P^*)$$
and
$\pi_{23}^*\P\tensor\pi_{12}^*\P^*=\pi_{23}^*\P^*\tensor\pi_{12}^*\P$.
But $\R\lhom(\hat E,\OO)=(-1)^*\R\FF(E^*_0)[n]$.
\end{proof}
\begin{cor}
$$\R\Phi(-)=(-1)^*E_0\Rconvol(-)[-n]\qquad\hbox{and}\qquad
\R\hat\Phi(-)=(-1)^*E_0^*\Rconvol(-),$$
where $\Rconvol$ denotes the (derived) convolution product:
$\R m_*(p_1^*(-)\tensor p_2^*(-))$ and
$m,p_1,p_2:\T\cross\T\to\T$ are the multiplication map and the projections
onto the factors respectively. In particular, $E^*_0\Rconvol E_0\cong\OO$.
\end{cor}
\begin{proof}
This follows immediately from the lemma and the convolution theorem for the
Mukai transform (see \refb{Muk2}{3.7}).
\end{proof}
\begin{proof} (of proposition)
Observe that the fibres of $\EE$ over $\{a\}\cross\{b\}$ are given by
$H^{n-i}(\R\FF(E_0)\tensor\P_{b}\tensor\P_{-a})$. This shows that
$\EE|_{\T\cross\{b\}}\cong E_b$ and so $\EE^{\oplus\sigma}$ is isomorphic
to the semi-universal bundle on $\T\cross\M$ via
$\T\mathbin{\mathop{\rightarrow}\limits^\sim}\M$. This
implies that $\EE$ is universal. But since $\EE_{\{a\}\cross\T}\cong
E_{-a}$ we see that it is also bi-universal.
\end{proof}
\begin{thm}\label{t:eulercl}
If $\M$ is a moduli space of stable bundles over an abelian variety
$\T$ of dimension $n$ satisfying WIT on which
$\T$ acts freely and effectively by pullback, then for each $E\in\M$ the
transform sheaf $\R\FF(E)$ is a line bundle. In particular, $|\chi(E)|=1$.
\end{thm}
\begin{proof}
The proposition shows that $\R\Phi$ and $\R\hat\Phi$ are of
Fourier-Mukai type. If we substitute the expressions of
\lemref{l:phiisf} into $\R\hat\Phi\compo\R\Phi=\Id$ then we obtain
$$\R\FF(E^*_0)\tensor\R\FF(E_0)=\OO.$$
This implies that $\R\FF(E_0)$ is an invertible sheaf.
\end{proof}

This allows us to generalise \thmref{t:modideal} to arbitrary
polarised abelian varieties and so we can state the following theorem.
\begin{thm}\label{t:hilbismod}
Let $\T$ be an abelian variety with a polarisation $\ell$. Then there
is a non-empty component $\M$ of the moduli space of simple sheaves on
$\T$ with Chern character 
$$\Bigl(\frac{m}{n!}(\hat\ell^n)^*-(-1)^n,-\frac{m}{(n-1)!}(\hat\ell^{n-1})^*,
\ldots,(-1)^{n-1}m\hat\ell^*,(-1)^nm1^*\Bigl),$$
where $\alpha\mapsto\alpha^*$ denotes $H^i(\dT,\C)\cong H^{2n-i}(\T,\C)$,
which is isomorphic to $\Hilb^m\T\cross\dT$. In particular, each Hilbert
scheme of points on an abelian variety arises as a moduli space of
simple sheaves on that variety. 
\end{thm}
\begin{proof} The proof is essentially the same as that of
\thmref{t:modideal}. The equivalent statement to \propref{p:7.1} holds because
$R^0\FF(L)$ has a natural Hermitian-Einstein connection via the flat
connection on $\PP(R^0\FF(L))$. Then the moduli space of Gieseker stable
vector bundles of Chern character $((\ell^n)^*/n!,(\ell^{n-1})^*/(n-1)!,\ldots,1^*)$
gives rise to a Fourier-Mukai transform $R\Phi$ as before. An ideal sheaf
$\I_S$ of a zero-dimensional subscheme is $\Phi$-WIT and the transform is a
simple sheaf by the Parseval Theorem.
\end{proof}
Analogous results also hold in the case of the K3 surface (see
\cite{BM}). This theorem strengthens the results of Mukai which give a series
of birational isomorphisms between Hilbert schemes of points and components of
the moduli of simple sheaves on abelian varieties (see
\cite{Muk5} Theorem 2.7, Theorem 2.17 and Theorem 2.20).
\section{Discussion}
The Fourier-Mukai transforms are very useful tools in the study of moduli
spaces of simple or stable sheaves as well as to the study of more direct
questions about the geometry of the underlying varieties. It is therefore an
important programme to identify them for any given variety with trivial
canonical bundle. As the application demonstrates it is possible to identify
many moduli spaces with Hilbert schemes of points. This lends weight to the
conjecture that all the components of the moduli spaces of simple sheaves on
projective varieties with trivial canonical bundle are isomorphic to
(deformations of) punctual Hilbert schemes.

The theory of Fourier-Mukai transforms presents a number of immediate
conjectures:
\begin{conjecture}
If $\R\Phi$ is a Fourier-Mukai functor from $X$ to $Y$ then $X$ is a
holomorphic deformation of $Y$.
\end{conjecture}
The Mukai transform shows that $X$ need not be naturally isomorphic to $Y$
when $X$ is an abelian variety with no principal polarisations. But the
geometry of $X$ and $Y$ are identical in this case. This is essentially the
question of the extent to which $D(X)$ determines $X$. For a K3 surface, no
examples are currently known for which $Y$ is not identical to $X$. For
abelian varieties, it is possible to find transforms when $X$ and $Y$ are
simply isogenous. 

A more specific question which arises based on the theorems above is
\begin{problem}
Given a Chern character, find a Fourier-Mukai transform for which each
simple sheaf in a component of the moduli space of simple sheaves with
this Chern character satisfies $\Phi$-WIT. This would be particularly
useful for abelian varieties where the sheaves might be half-WITs.
\end{problem}
It would also be important to know when strongly bi-universal sheaves exist in
more general contexts because the resulting transforms are still of Verdier
type and may even have one sided inverses.

\begin{conjecture}
Given a smooth projective variety $X$ with $K_X=\OO_X$ and $Y$ a
projective component of $\Spl(X)$ of dimension $\dim X$ with a
universal sheaf $\EE$ on $X\cross Y$. Then $\EE$ is
bi-universal.
\end{conjecture}
\noindent All currently known examples satisfy this conjecture, for example,
this is always true for an abelian surface.

Another interesting problem is to construct such Fourier-Mukai
transforms for Calabi-Yau three-folds. One could then study the
analytic versions of the transforms (in analogy with the Nahm
transform for instantons on complex tori) and use these to solve the
Hermitian-Yang-Mills equations on such 3-folds. An obvious question then
arises about whether one can find a Fourier-Mukai transform for which $Y$ is
a mirror of $X$.

\end{document}